\documentclass[12pt]{article}
\usepackage{amssymb,amsmath,amsthm,amsfonts,enumerate, bbm, mathdots}

\usepackage{latexsym}
\usepackage{amscd}
\usepackage{stmaryrd}
\usepackage{color}
\usepackage[all]{xypic}
\usepackage{epsfig}
\usepackage{graphics}
\usepackage{ifthen}
\usepackage{varioref}
\usepackage{rotating}
\usepackage{extarrows}
\usepackage{cite}
\usepackage{mathrsfs}

\numberwithin{equation}{section}

\theoremstyle{plain}
\newtheorem{thm}{Theorem}[section]
\newtheorem{cor}[thm]{Corollary}
\newtheorem{lem}[thm]{Lemma}
\newtheorem{prop}[thm]{Proposition}
\newtheorem{rem}[thm]{Remark}
\newtheorem{Def}[thm]{Definition}

\makeatletter
\newenvironment{proofof}[1]{\par
  \pushQED{\qed}%
  \normalfont \topsep6\p@\@plus6\p@\relax
  \trivlist
  \item[\hskip\labelsep
        \bfseries
    Proof of #1\@addpunct{.}]\ignorespaces
}{%
  \popQED\endtrivlist\@endpefalse
}
\makeatother

\definecolor{darkgreen}{rgb}{0.0625,0.64,0.0625}
\usepackage[
            pdfstartview=FitH,
            CJKbookmarks=true,
            bookmarksnumbered=true,
            bookmarksopen=true,
            colorlinks,
            pdfborder=001,
            linkcolor=blue,
            anchorcolor=green,
            citecolor=red
            ]{hyperref}

\newfont{\scyr}{wncyr10 scaled 550}

\def\proof{\noindent {\bf Proof.\;}}

\def\wt{\operatorname{wt}}

\def\sgn{\operatorname{sgn}}
\def\max{\operatorname{max}}

\allowdisplaybreaks

\begin{document}
\title{Integrality and some evaluations of odd multiple harmonic sums\thanks{This work was supported by the Fundamental Research Funds for the Central Universities (grant number 22120210552).}}
\date{\small ~ \qquad\qquad School of Mathematical Sciences, Tongji University \newline No. 1239 Siping Road,
Shanghai 200092, China}

\author{Zhonghua Li\thanks{E-mail address: zhonghua\_li@tongji.edu.cn} ~and ~Zhenlu Wang\thanks{E-mail address: zlw@tongji.edu.cn}}
\maketitle
\begin{abstract}
In 2015, S. Hong and C. Wang proved that none of the elementary symmetric functions of $1,1/3,\ldots,1/(2n-1)$ is an integer when $n\geq 2$. In 2017, Kh. Pilehrood, T. Pilehrood and R. Tauraso proved that the multiple harmonic sums $H_n(s_1,\ldots,s_r)$ are never integers with exceptions of $H_1(s_1)=1$ and $H_3(1,1)=1$. They also proved that the multiple harmonic star sums are never integers when $n\geq 2$. In this paper, we consider the odd multiple harmonic sums and the odd multiple harmonic star sums and show that none of these sums is an integer with exception of the trivial case. Besides, we give evaluations of the odd (alternating) multiple harmonic sums with depth one.
\end{abstract}

{\small

}

\section{Introduction}
Let $\mathbb{N}$ be the set of positive integers. For any $r,n\in \mathbb{N}$ with $r\leq n$ and $\mathbf{s}=(s_1,s_2,\ldots,s_r)\in\mathbb{N}^r$, the multiple harmonic sum and the multiple harmonic star sum are defined respectively by
\begin{align}\label{MHS Def}
H_n(\mathbf{s})=H_n(s_1,s_2,\ldots,s_r):=\sum\limits_{1\leq k_1<k_2<\cdots<k_r\leq n}\prod_{j=1}^{r}\frac{1}{k_j^{s_j}},
\end{align}
and
\begin{align}\label{MsHS Def}
H_n^{\star}(\mathbf{s})=H_n^{\star}(s_1,s_2,\ldots,s_r):=\sum\limits_{1\leq k_1\leq k_2\leq\cdots\leq k_r\leq n}\prod_{j=1}^{r}\frac{1}{k_j^{s_j}}.
\end{align}
The odd multiple harmonic sum and the odd multiple harmonic star sum are defined respectively by
\begin{align}\label{OddMHS Def}
\overline{H}_n(\mathbf{s})=\overline{H}_n(s_1,s_2,\ldots,s_r):=\sum\limits_{0\leq k_1<k_2<\cdots<k_r\leq n-1}\prod_{j=1}^{r}\frac{1}{(2k_j+1)^{s_j}},
\end{align}
and
\begin{align}\label{OddMsHS Def}
\overline{H}_n^{\star}(\mathbf{s})=\overline{H}_n^{\star}(s_1,s_2,\ldots,s_r):=\sum\limits_{0\leq k_1\leq k_2\leq \cdots\leq k_r\leq n-1}\prod_{j=1}^{r}\frac{1}{(2k_j+1)^{s_j}}.
\end{align}
We call $r$ the depth and $\wt(\mathbf{s}):=\sum\limits_{j=1}^rs_j$ the weight. Besides, we denote the sequence of $a$ with $r$ repetitions as $\left(\{a\}^r\right)$.

The multiple harmonic sums \eqref{MHS Def} and \eqref{OddMHS Def} are of certain interest because by taking the limit as $n$ goes to $\infty$ when $s_r>1$, we get
\begin{align*}
\lim\limits_{n\rightarrow\infty}H_n(s_1,\ldots,s_r)=\zeta(s_1,\ldots,s_r)\quad\text{and}\quad \lim\limits_{n\rightarrow\infty}\overline{H}_n(s_1,\ldots,s_r)=t(s_1,\ldots,s_r),
\end{align*}
which are the so-called multiple zeta value \cite{Hoffman92,Zagier} and Hoffman's multiple $t$-value \cite{Hoffman}, respectively. Also by taking the limit, the star versions \eqref{MsHS Def} and \eqref{OddMsHS Def} deduce the multiple zeta-star value \cite{Hoffman92} and Hoffman's multiple $t$-star value  \cite{Hoffman} respectively.

A well-known result in elementary number theory says that for any positive integers $n>1$,  the harmonic sum $H_n(1)$ is not an integer. In 1946, P. Erd{\H o}s and I. Niven \cite{Erdos-Niven} proved that the nested multiple harmonic sums $H_n(\{1\}^r)$ can take integer values only for a finite number of positive integers $n$, and they also mentioned that there are only finitely many integers $n$ for which one or more of the elementary symmetric functions of $1/a,1/(a+b),\ldots,1/(a+bn)$ are integers, where $a$ and $b$ are given positive integers. Y. Chen and M. Tang \cite{Chen-Tang} proved that $H_n(\{1\}^r)$ is an integer only for $(n,r)=(1,1)$ and $(n,r)=(3,2)$. Kh. Pilehrood, T. Pilehrood and R. Tauraso \cite{P-P-T} got a more general result which says that for $s_1,\ldots,s_r\in\mathbb{N}$, $H_n(s_1,\ldots,s_r)$ is never an integer with exceptions of $H_1(s_1)=1$ and $H_3(1,1)=1$, and they also proved that $H_n^{\star}(s_1,\ldots,s_r)$ is never an integer with exception of $H_1^{\star}(s_1)=1$. For the odd multiple harmonic sums, S. Hong and C. Wang \cite{Hong-Wang-2015} proved that $\overline{H}_n(\{1\}^r)$ is not an integer for any integer $n\geq 2$. For other works such as the binomial identities and the congruences  involving the (odd) multiple harmonic sums, one can refer to \cite{P-P-T-Def} or the book \cite{Zhao} of J. Zhao and the references therein.

In this paper, using similar method as in \cite{P-P-T}, we prove the following results indicates the integrality of the odd multiple harmonic star sums and the odd multiple harmonic sums.

\begin{thm}\label{star-Theorem integrity}
Let $r,n,s_1,\ldots,s_r\in \mathbb{N}$ such that $r\leq n$. Then $\overline{H}_n^{\star}(s_1,\ldots,s_r)$ is never an integer when $n\geq 2$.
\end{thm}

\begin{thm}\label{Theorem integrity}
Let $r,n,s_1,\ldots,s_r\in \mathbb{N}$ such that $r\leq n$. Then $\overline{H}_n(s_1,\ldots,s_r)$ is never an integer when $n\geq 2$.
\end{thm}

More generally, one can define the odd alternating multiple harmonic sums. Let $\mathbb{D}=\mathbb{N}\cup\overline{\mathbb{N}}$, where $\overline{\mathbb{N}}=\{\overline{s}\mid s\in \mathbb{N}\}$. Define the absolute value function $|\cdot|$ on $\mathbb{D}$ by $|s|=|\overline{s}|=s$ for all $s\in\mathbb{N}$ and the sign function by $\sgn(s)=1$ and $\sgn(\overline{s})=-1$ for all $s\in \mathbb{N}$. Let $n,r\in \mathbb{N}$ with $r\leq n $ and $\mathbf{s}=(s_1,s_2,\ldots,s_r)\in\mathbb{D}^r$, the odd alternating multiple harmonic sum $\overline{H}_n(\mathbf{s})$ is defined by
$$\overline{H}_n(\mathbf{s})=\overline{H}_n(s_1,s_2,\ldots,s_r):=\sum\limits_{0\leq k_1<k_2<\cdots<k_r\leq n-1}\prod_{j=1}^{r}\frac{\sgn(s_j)^{k_j}}{(2k_j+1)^{|s_j|}}.$$

A formula known to Euler says that for any $n\in\mathbb{N}$, it holds
$$H_n(1)=\sum\limits_{k=1}^n\frac{(-1)^{k-1}}{k}\binom{n}{k}.$$
In this paper, we generalize the above binomial identity and give the evaluations of $\overline{H}_n(s)$ and $\overline{H}_n(\overline{s})$ for any $n,s\in\mathbb{N}$. In the evaluation formulas, the special values of the generalized hypergeometric functions appear.

The paper is organized as follows. In Section \ref{proof}, we first give the proof of Theorem \ref{star-Theorem integrity}. Then we show several lemmas which are needed for the proof of Theorem \ref{Theorem integrity} and finally complete the proof. In Section \ref{values}, we give the evaluations of $\overline{H}_n(s)$ and $\overline{H}_n(\overline{s})$ for any $n,s\in\mathbb{N}$.

In this paper, all the numerical computations were performed by using Matlab.

\section{Proof of Theorem \ref{star-Theorem integrity} and Theorem \ref{Theorem integrity}}\label{proof}

For a prime $p$, let $v_p(a)$ be the $p$-adic valuation of the rational number $a$, that is, if $a=p^nq_1/q_2$ with $n,q_1,q_2\in\mathbb{Z}$ and $q_1,q_2$ are coprime with $p$, then $v_p(a)=n$. It is a fact that the $p$-adic valuation satisfies the inequality $$v_p(a+b)\geq \min(v_p(a),v_p(b)),$$ where the equality holds for $v_p(a)\neq v_p(b)$.

The integrality of the odd multiple harmonic star sums is simple to prove.

\begin{proofof}{Theorem \ref{star-Theorem integrity}}
Let $\mathbf{s}=(s_1,s_2,\ldots,s_r)$. Since $n\geq 2$, by Bertrand's postulate, there exists a prime $p$ such that $n<p<2n$. Hence, we obtain that
\begin{align*}
\overline{H}_n^{\star}(\mathbf{s})&=\sum\limits_{0\leq k_1\leq k_2\leq \cdots\leq k_r\leq n-1\atop \exists j,p\neq 2k_j+1}\prod_{j=1}^{r}\frac{1}{(2k_j+1)^{s_j}}+\sum\limits_{0\leq k_1\leq k_2\leq \cdots\leq k_r\leq n-1\atop \forall j,p=2k_j+1}\prod_{j=1}^{r}\frac{1}{(2k_j+1)^{s_j}}\\
&=\frac{a}{bp^t}+\frac{1}{p^{\wt(\mathbf{s})}},
\end{align*}
where $a,b$ are positive integers with $(b,p)=1$, $t$ is a nonnegative integer with $t<\wt(\mathbf{s})$. Since $v_p\left(\frac{a}{bp^t}\right)>v_p\left(\frac{1}{p^{\wt(\mathbf{s})}}\right)$, we have $v_p\left(\overline{H}_n^{\star}(\mathbf{s})\right)=-\wt(\mathbf{s})<0$, which implies $\overline{H}_n^{\star}(\mathbf{s})$ is not an integer.
\end{proofof}


The proof of Theorem \ref{Theorem integrity} is more complicated. We need several lemmas. In order to compare values of odd multiple harmonic sums of the same depth, the following definition and lemma will be useful.

\begin{Def}[{\cite[Definition 1]{P-P-T}}]
Let $\mathbf{s}=(s_1,\ldots,s_r)$ and $\mathbf{t}=(t_1,\ldots,t_r)$ be two $r$-tuples of positive integers. If $\wt(\mathbf{s})\geq \wt(\mathbf{t})$, and $s_1\leq t_1,\ldots,s_l\leq t_l,s_{l+1}\geq t_{l+1},\ldots,s_r\geq t_r$ for some $0\leq l\leq r-1$, we say that $\mathbf{s}\geq \mathbf{t}$.
\end{Def}

\begin{lem}\label{order of s,t}
For $r,n\in \mathbb{N}$ such that $r\leq n$, if $\mathbf{s}=(s_1,\ldots,s_r)$ and $\mathbf{t}=(t_1,\ldots,t_r)$ are two $r$-tuples of positive integers such that $\mathbf{s}\geq \mathbf{t}$, then it holds $$\overline{H}_n(s_1,\ldots,s_r)\leq \overline{H}_n(t_1,\ldots,t_r).$$
\end{lem}

\proof
Let $0\leq l\leq r-1$ such that $s_1\leq t_1,\ldots,s_l\leq t_l,s_{l+1}\geq t_{l+1},\ldots,s_r\geq t_r$. For fixed $0\leq k_1<\cdots<k_r\leq n-1$, we have
\begin{align*}
(2k_1+1)^{t_1-s_1}\cdots(2k_l+1)^{t_l-s_l}&\leq (2k_l+1)^{t_1+\cdots+t_l-(s_1+\cdots+s_l)}\\
&=(2k_l+1)^{t_1+\cdots+t_r-(s_1+\cdots+s_r)}\cdot(2k_l+1)^{s_{l+1}-t_{l+1}+\cdots+s_r-t_r}\\
&\leq (2k_l+1)^{s_{l+1}-t_{l+1}+\cdots+s_r-t_r}\\
&\leq (2k_{l+1}+1)^{s_{l+1}-t_{l+1}}\cdots(2k_r+1)^{s_r-t_r}.
\end{align*}
Therefore, we get
$$\frac{1}{(2k_1+1)^{s_1}\cdots(2k_r+1)^{s_r}}\leq \frac{1}{(2k_1+1)^{t_1}\cdots(2k_r+1)^{t_r}},$$
which implies the desired result.\qed

We then show that the odd multiple harmonic sum is not an integer if the depth is large enough.

\begin{lem}\label{Lemma 2<=r<}
Let $r,n\in \mathbb{N}$ such that
$$e\left(\frac{1}{2}\log(2n-1)+1\right)\leq r\leq n.$$
Then for any $s_1,\ldots,s_r\in\mathbb{N}$, $\overline{H}_n(s_1,\ldots,s_r)$ is not an integer.
\end{lem}

\proof
Notice that
$$\overline{H}_n(1)=\sum\limits_{i=0}^{n-1}\frac{1}{2i+1}<1+\int_0^{n-1}\frac{1}{2x+1}dx=\frac{1}{2}\log(2n-1)+1.$$
By using Lemma \ref{order of s,t}, we get
\begin{align}
\overline{H}_n(s_1,\ldots,s_r)\leq \overline{H}_n(\{1\}^r)\leq \frac{\left(\overline{H}_n(1)\right)^r}{r!}\leq \frac{\left(\frac{1}{2}\log(2n-1)+1\right)^r}{r!}.
\end{align}
Since $e\left(\frac{1}{2}\log(2n-1)+1\right)\leq r$ and $\frac{r^r}{r!}<e^r$, we have $\overline{H}_n(s_1,\ldots,s_r)<1$, which implies the result.
\qed

The following lemma is useful in the proof of Theorem \ref{Theorem integrity}.

\begin{lem}\label{p condition}
Let $r,n\in \mathbb{N}$ such that $r\leq n$. Suppose that there exists a prime $p>r+1$ satisfying that
$$\frac{2n}{r+1}\leq p<\frac{2n}{r},$$
then for any $s_1,\ldots,s_r\in\mathbb{N}$, $\overline{H}_n(s_1,\ldots,s_r)$ is not an integer.
\end{lem}

\proof
Since $0<p<\cdots<rp<2n\leq (r+1)p<p^2$ and $s_1,\ldots,s_r$ are positive integers, we have
\begin{align*}
\overline{H}_n(s_1,\ldots,s_r)&=\sum\limits_{0\leq k_1<\cdots<k_r\leq n-1\atop \exists j, p\nmid 2k_j+1}\prod_{j=1}^r\frac{1}{(2k_j+1)^{s_j}}+\sum\limits_{0\leq k_1<\cdots<k_r\leq n-1\atop \forall j, p\mid 2k_j+1}\prod_{j=1}^r\frac{1}{(2k_j+1)^{s_j}}\\
&=\frac{a}{bp^t}+\frac{1}{cp^{\wt(\mathbf{s})}},
\end{align*}
where $\mathbf{s}=(s_1,\ldots,s_r)$, $a,b,c$ are positive integers with $(b,p)=(c,p)=1$, $t$ is a nonnegative integer with $t<\wt(\mathbf{s})$.
Therefore, $v_p\left(\overline{H}_n(s_1,\ldots,s_r)\right)=-\wt(\mathbf{s})<0$, which implies $\overline{H}_n(s_1,\ldots,s_r)$ is not an integer.
\qed

Now we show that in the case of $1\leq r<e\left(\frac{1}{2}\log(2n-1)+1\right)$, $\overline{H}_n(s_1,\ldots,s_r)$ is not an integer if $n$ is large enough. We need a lemma.

\begin{lem}[{\cite[Theorem 1.9]{Dusart}}]\label{P. Dusart theorem}
For any real numbers $x\geq 3275$, there exists a prime $p$ satisfying $x<p\leq x\left(1+\frac{1}{2\log^2x}\right)$.
\end{lem}

Using the above lemmas, we have the following result.

\begin{lem}\label{Lemma n>=30000}
Let $r,n\in \mathbb{N}$ satisfying $n\geq 30000$ and $r<e\left(\frac{1}{2}\log(2n-1)+1\right)$. Then for any $s_1,\ldots,s_r\in \mathbb{N}$, we have $\overline{H}_n(s_1,\ldots,s_r)$ is not an integer.
\end{lem}

\proof
Let $x=\frac{2n}{r+\frac{1}{2}}$. By Lemma \ref{P. Dusart theorem}, if
\begin{align}\label{x>=3275}
x\geq 3275,
\end{align}
and
\begin{align}\label{log(x)>r}
\log^2x>r,
\end{align}
then there exists a prime $p$ satisfying $\frac{2n}{r+\frac{1}{2}}<p<\frac{2n}{r}$, which implies $\frac{2n}{r+1}<p<\frac{2n}{r}$. If $p$ satisfies
\begin{align}\label{p>r+1}
p>r+1,
\end{align}
then $\overline{H}_n(s_1,\ldots,s_r)$ is not an integer by Lemma \ref{p condition}. Hence we only need to show \eqref{x>=3275}, \eqref{log(x)>r} and \eqref{p>r+1}.

First, we prove that $x\geq 3275$. Let
$$p(x)=2x-3275\left(\frac{e}{2}\log(2x-1)+e+\frac{1}{2}\right).$$
As
\begin{align*}
p'(x)=2-\frac{3275e}{2x-1}>0
\end{align*}
for all $x\geq 30000$ and $p(30000)>0$, we get $p(x)>0$ for all $x\geq 30000$. Since $r<e\left(\frac{1}{2}\log(2n-1)+1\right)$, we prove \eqref{x>=3275}.

Then, we prove that $\log^2x>r$. It is sufficient to prove that
\begin{align}\label{log(x),r prove}
\log^22n-2\log2n\cdot\log\left(r+\frac{1}{2}\right)>r.
\end{align}\label{log(x),r prove}
Since $r<e\left(\frac{1}{2}\log(2n-1)+1\right)$, it is enough to prove that
$$\log2n-2\log\left(\frac{e}{2}\log2n+e+\frac{1}{2}\right)>\frac{e}{2}+\frac{e}{\log2n}.$$
Let $$q(x)=x-2\log\left(\frac{e}{2}x+e+\frac{1}{2}\right)-2.$$
As
\begin{align*}
q'(x)=1-\frac{e}{\frac{e}{2}x+e+\frac{1}{2}}>0
\end{align*}
for all $x>0$ and $q(\log(2\cdot30000))>0$, we obtain $q(x)>0$ for all $x\geq \log(2\cdot30000)$. Since for all $n\geq 30000$, $$\frac{e}{2}+\frac{e}{\log2n}\leq \frac{e}{2}+\frac{e}{\log(2\cdot30000)}<2,$$
\eqref{log(x)>r} is proved.

Finally, we prove that $p>r+1$. Since $p>\frac{2n}{r+1}$, it is sufficient to prove that  $\sqrt{2n}>r+1$. Let
$$t(x)=\sqrt{2x}-\frac{e}{2}\log(2x-1)-e-1.$$
Since
$$t'(x)=\frac{2x-1-e\sqrt{2x}}{\sqrt{2x}(2x-1)}>0$$
for all $x\geq 30000$ and $t(30000)>0$, we have $t(x)>0$ for all $x\geq 30000$. Therefore, \eqref{p>r+1} is proved.
\qed

For $r\in \mathbb{N}$, set
$$B_r=\bigcup\limits_{p\in\mathbb{P}}(rp,(r+1)p],$$
where $\mathbb{P}$ denotes the set of all primes. As a variant of \cite[Lemma 1]{P-P-T}, we can prove that $B_r$ is cofinite in $\mathbb{N}$ similarly. Also, we can prove that
$$\max(\mathbb{N}\setminus B_r)=\max(\mathbb{N}\setminus A_r)+1,$$
where $A_r=\bigcup\limits_{p\in\mathbb{P}}[rp,(r+1)p)$ was defined in  \cite{P-P-T}. Let $n_r=[\frac{\max(\mathbb{N}\setminus B_r)}{2}]+1$. Then using the table given in \cite[Lemma 1]{P-P-T}, we have the following table for some values of $n_r$:

\vspace{0.45cm}
\begin{center}
\begin{tabular}{|c|c|c|c|c|c|c|c|c|c|c|}
\hline
$r$&1&2&3&4&5&6&7&8&9&10\\
\hline
$n_r$&2&12&17&59&73&112&130&213&572&636\\
\hline
$r$&11&12&13&14&15&16&17&18&19&20\\
\hline
$n_r$&699&763&826&1044&1118&1193&1794&2008&2119&2231\\
\hline
\end{tabular}
\end{center}
\vspace{0.45cm}\

\begin{rem}\label{Remark-N-r}
If $r\geq 1$ and $n\geq n_r$, then $2n\in B_r$. Hence there is a prime $p$, such that $\frac{2n}{r+1}\leq p<\frac{2n}{r}$. For $2\leq r\leq 20$, we may verify that $2n_r>(r+1)^2$, which implies $p>r+1$. Thus using Lemma \ref{p condition}, we find that $\overline{H}_n(s_1,\ldots,s_r)$ is not an integer for $2\leq r\leq 20$, $n\geq n_r$ and $s_1,\ldots,s_r\in\mathbb{N}$.
\end{rem}

Similarly as in \cite[Lemma 3]{P-P-T}, we have to find an upper bound for the index $s_1$.

\begin{lem}\label{Lemma M}
Given $2\leq r\leq n$ and $(s_2,\ldots,s_r)\in \mathbb{N}^{r-1}$, there exists an integer $N$ (which depends on $(s_2,\ldots,s_r)$ and $n$) such that $\overline{H}_n(s_1,\ldots,s_r)$ is not an integer for any positive integer $s_1>N$.
\end{lem}

\proof
Let $(s_1,\ldots,s_r)\in \mathbb{N}^{r}$. If $r=n$, it is obvious that $\overline{H}_n(s_1,\ldots,s_r)$ is not an integer for all $s_1\geq 1$. If $r<n$, we set
\begin{align*}
\overline{H}_n(s_1,\ldots,s_r)=\sum\limits_{k_1=0}^{n-r}\frac{c_{k_1}}{(2k_1+1)^{s_1}}\quad\text{where}\quad c_{k_1}=\sum\limits_{k_1<k_2<\cdots<k_r\leq n-1}\prod\limits_{j=2}^r\frac{1}{(2k_j+1)^{s_j}}.
\end{align*}
By Bertrand's postulate, there exists at least one prime in $(n-r+1,2n-2r+2)$. Let $p$ be such largest prime, then $2p>2n-2r+2$. We get
\begin{align*}
\overline{H}_n(s_1,\ldots,s_r)=\sum\limits_{k_1=0\atop p=2k_1+1}^{n-r}\frac{c_{k_1}}{(2k_1+1)^{s_1}}+\sum\limits_{k_1=0\atop p\neq 2k_1+1}^{n-r}\frac{c_{k_1}}{(2k_1+1)^{s_1}}
&=\frac{c_{\frac{p-1}{2}}}{p^{s_1}}+\sum\limits_{k_1=0\atop p\neq 2k_1+1}^{n-r}\frac{c_{k_1}}{(2k_1+1)^{s_1}}.
\end{align*}
Define
$$N(n,s_2,\ldots,s_r):=\max\left(v_p(c_{\frac{p-1}{2}}),v_p(c_{\frac{p-1}{2}})-\min\limits_{0\leq k_1\leq n-r\atop p\neq 2k_1+1}v_p(c_{k_1})\right).$$
Assume that $s_1>N(n,s_2,\ldots,s_r)$, then we have
\begin{align}\label{s_1 1}
s_1>v_p(c_{\frac{p-1}{2}})
\end{align}
and
\begin{align}\label{s_1 2}
s_1>v_p(c_{\frac{p-1}{2}})-\min\limits_{0\leq k_1\leq n-r\atop p\neq 2k_1+1}v_p(c_{k_1}).
\end{align}
By \eqref{s_1 2}, we have
\begin{align*}
v_p\left(\sum\limits_{k_1=0\atop p\neq 2k_1+1}^{n-r}\frac{c_{k_1}}{(2k_1+1)^{s_1}}\right)\geq \min\limits_{0\leq k_1\leq n-r\atop p\neq 2k_1+1}v_p\left(\frac{c_{k_1}}{(2k_1+1)^{s_1}}\right)=\min\limits_{0\leq k_1\leq n-r\atop p\neq 2k_1+1}v_p(c_{k_1})>v_p\left(\frac{c_{\frac{p-1}{2}}}{p^{s_1}}\right).
\end{align*}
Hence, by \eqref{s_1 1}, we get
$$v_p\left(\overline{H}_n(s_1,\ldots,s_r)\right)=\min\left(v_p\left(\frac{c_{\frac{p-1}{2}}}{p^{s_1}}\right),v_p\left(\sum\limits_{k_1=0\atop p\neq 2k_1+1}^{n-r}\frac{c_{k_1}}{(2k_1+1)^{s_1}}\right)\right)=v_p\left(\frac{c_{\frac{p-1}{2}}}{p^{s_1}}\right)<0,$$
which implies $\overline{H}_n(s_1,\ldots,s_r)$ is not an integer.
\qed

Now we come to prove Theorem \ref{Theorem integrity}.

\begin{proofof}{Theorem \ref{Theorem integrity}}
Let $\mathbf{s}=(s_1,\ldots,s_r)$. The case of $r=1$ follows from Theorem \ref{star-Theorem integrity}. Then by Lemma \ref{Lemma 2<=r<} and Lemma \ref{Lemma n>=30000}, we may assume that $2\leq r<e\left(\frac{1}{2}\log(2n-1)+1\right)$ and $n<30000$. In this case, sine $e\left(\frac{1}{2}\log(2n-1)+1\right)<18$ and Remark \ref{Remark-N-r}, it is enough to consider the integrality of $\overline{H}_n(s_1,\ldots,s_r)$ for $2\leq r\leq 17$ and $n<n_r$.

Let $11\leq r\leq 17$. For $n<n_r$, we have
\begin{align*}
\overline{H}_n(s_1,\ldots,s_r)<\overline{H}_{n_r}(s_1,\ldots,s_r)\leq\overline{H}_{n_r}(\{1\}^r)\leq\frac{\left(\frac{1}{2}\log(2n_r-1)+1\right)^r}{r!}<1,
\end{align*}
where the last inequality can be verified numerically. Therefore, $\overline{H}_n(s_1,\ldots,s_r)$ is not an integer.

Let $5\leq r\leq 10$. We have the following estimations
$$\overline{H}_{n_{10}}(s_1,\ldots,s_{10})\leq\overline{H}_{n_{10}}(\{1\}^{10})<0.06243<1,$$
$$\overline{H}_{n_9}(s_1,\ldots,s_9)\leq\overline{H}_{n_9}(\{1\}^{9})<0.17796<1,$$
$$\overline{H}_{n_8}(s_1,\ldots,s_8)\leq\overline{H}_{n_8}(\{1\}^{8})<0.12835<1,$$
$$\overline{H}_{n_7}(s_1,\ldots,s_7)\leq\overline{H}_{n_7}(\{1\}^{7})<0.20280<1,$$
$$\overline{H}_{n_6}(s_1,\ldots,s_6)\leq\overline{H}_{n_6}(\{1\}^{6})<0.51825<1,$$
$$\overline{H}_{n_5}(s_1,\ldots,s_5)\leq\overline{H}_{n_5}(\{1\}^{5})<0.85442<1.$$
Hence, for $n<n_r$, $\overline{H}_n(s_1,\ldots,s_r)$ is not an integer.

Let $r=2,2\leq n<n_2=12$. We verify numerically that $\overline{H}_{n}(1,1)$ is not an integer. Since $\overline{H}_{12}(1,2)<0.27273<1$, we get $$\overline{H}_{n}(1,s_2)\leq\overline{H}_{n}(1,2)<\overline{H}_{12}(1,2)<1$$
for all $s_2>1$, which implies $\overline{H}_{n}(1,s_2)$ is not an integer if $s_2>1$. We can verify that the upper bound $N(n,1)=1$ for all $3\leq n\leq 11$ in the proof of Lemma \ref{Lemma M}. Hence $\overline{H}_{n}(s_1,1)$ is not an integer for $s_1>1$ and $3\leq n\leq 11$ by Lemma \ref{Lemma M}. Therefore we proved that $\overline{H}_{n}(s_1,s_2)$ is never an integer for all $s_1,s_2\in\mathbb{N}$.

Let $r=3,3\leq n<n_3=17$. We verify numerically that $\overline{H}_{n}(1,1,1)$ is not an integer. Since $\overline{H}_{17}(1,2,1)<0.22216<1$, we get
\begin{align*}
\overline{H}_{n}(s_1,s_2,s_3)\leq\overline{H}_{n}(1,2,1)<\overline{H}_{17}(1,2,1)<1
\end{align*}
for all $s_1,s_3\geq 1$ and $s_2>1$. Since $\overline{H}_{17}(1,1,2)<0.08552<1$, we have
\begin{align*}
\overline{H}_{n}(s_1,s_2,s_3)\leq\overline{H}_{n}(1,1,2)<\overline{H}_{17}(1,1,2)<1,
\end{align*}
for all $s_1,s_2\geq 1$ and $s_3>1$. As we can verify $\max\limits_{4\leq n\leq 16}N(n,1,1)=1$, we find from Lemma \ref{Lemma M} that $\overline{H}_{n}(s_1,1,1)$ is not an integer for $s_1>1$. Thus, we proved that $\overline{H}_{n}(s_1,s_2,s_3)$ is never an integer for all $s_1,s_2,s_3\in\mathbb{N}$.

Let $r=4,4\leq n<n_4=59$. We verify numerically that $\overline{H}_{n}(1,1,1,1)$ is not an integer. Since $\overline{H}_{59}(1,2,1,1)<0.28433<1$, we get
\begin{align*}
\overline{H}_{n}(s_1,s_2,s_3,s_4)\leq\overline{H}_{n}(1,2,1,1)<\overline{H}_{59}(1,2,1,1)<1
\end{align*}
for all $s_1,s_3,s_4\geq1$ and $s_2>1$. Since $\overline{H}_{59}(1,1,2,1)<0.10452<1$, we get
\begin{align*}
\overline{H}_{n}(s_1,s_2,s_3,s_4)\leq\overline{H}_{n}(1,1,2,1)<\overline{H}_{59}(1,1,2,1)<1
\end{align*}
for all $s_1,s_2,s_4\geq1$ and $s_3>1$. Since $\overline{H}_{59}(1,1,1,2)<0.04060<1$, we get
\begin{align*}
\overline{H}_{n}(s_1,s_2,s_3,s_4)\leq\overline{H}_{n}(1,1,1,2)<\overline{H}_{59}(1,1,1,2)<1
\end{align*}
for all $s_1,s_2,s_3\geq1$ and $s_4>1$. As we have $\max\limits_{5\leq n\leq 58}N(n,1,1,1)=1$, we get that $\overline{H}_{n}(s_1,1,1,1)$ is not an integer for $s_1>1$ from Lemma \ref{Lemma M}. Therefore, we have $\overline{H}_{n}(s_1,s_2,s_3,s_4)$ is never an integer for all $s_1,s_2,s_3,s_4\in\mathbb{N}$.
\end{proofof}

\section{Some evaluations}\label{values}

In this section, we give evaluations of $\overline{H}_n(s)$  and $\overline{H}_n(\overline{s})$ for any $n,s\in\mathbb{N}$. The generalized hypergeometric functions are needed in the following discussion. For $a_1,\ldots,a_{s+1}$, $b_1,\ldots,b_s\in\mathbb{C}$ with none of  $b_i$ is zero or a negative integer, the generalized hypergeometric function is defined by
\begin{align*}
_{s+1}F_{s}(a_1,\ldots,a_{s+1};b_1,\ldots,b_s;x):=\sum\limits_{i=0}^\infty\frac{(a_1)_i\cdots(a_{s+1})_i}{(b_1)_i\cdots(b_s)_i}\frac{x^i}{i!},
\end{align*}
where
$$(a)_i:=\begin{cases}
1, &  i=0, \\
a(a+1)\cdots(a+i-1), & i>0.
\end{cases}$$
The above series is absolutely and uniformly convergent if $|x|<1$ and the convergence also extends over the unit circle if $\Re\left(\sum\limits_{i=1}^sb_i-\sum\limits_{i=1}^{s+1}a_i\right)>0$.

We have the following result.

\begin{lem}\label{lemma 3-1}
For any $n,s\in \mathbb{N}$, we have
\begin{align}\label{3-1}
\sum\limits_{k=0}^{n-1}\frac{x^{2k}}{(2k+1)^s}=\sum\limits_{k=1}^n(-1)^{k-1}\binom{n}{k}{_{s+1}F_s}\left(\left\{\frac{1}{2}\right\}^{s},1-k;\left\{\frac{3}{2}\right\}^{s};x^2\right).
\end{align}
\end{lem}

\proof
We prove the lemma by induction on $s$. Let
$$f_s(x)=\sum\limits_{k=0}^{n-1}\frac{x^{2k}}{(2k+1)^s}.$$
If $s=1$, then
$$f_1(x)=\sum\limits_{k=0}^{n-1}\frac{x^{2k}}{2k+1},$$
and we get that
$$(xf_1(x))'=\sum\limits_{k=0}^{n-1}x^{2k}=\sum\limits_{k=1}^n(-1)^{k-1}\binom{n}{k}(1-x^2)^{k-1}.$$
Integrating both sides from $0$ to $x$, we find the left-hand side becomes $xf_1(x)$ and the right-hand side turns to be
\begin{align*}
\sum\limits_{k=1}^n(-1)^{k-1}\binom{n}{k}\int_0^x(1-t^2)^{k-1}dt&=\sum\limits_{k=1}^n(-1)^{k-1}\binom{n}{k}\int_0^x\sum\limits_{l=0}^{k-1}(-1)^l\binom{k-1}{l}t^{2l}dt\\
&=\sum\limits_{k=1}^n(-1)^{k-1}\binom{n}{k}\sum\limits_{l=0}^{k-1}\frac{(-1)^l(k-1)\cdots(k-l)x^{2l+1}}{(2l+1)l!}\\
&=\sum\limits_{k=1}^n(-1)^{k-1}\binom{n}{k}\sum\limits_{l=0}^{k-1}\frac{\left(\frac{1}{2}\right)_l\cdot(1-k)_l\cdot x^{2l+1}}{\left(\frac{3}{2}\right)_l\cdot l!}\\
&=x\sum\limits_{k=1}^n(-1)^{k-1}\binom{n}{k}{_2}F_1\left(\frac{1}{2},1-k;\frac{3}{2};x^2\right).
\end{align*}
Thus, we prove the case of $s=1$. Now Assume that \eqref{3-1} holds for $s$. Integrating both sides of \eqref{3-1} from $0$ to $x$, we find the left-hand side becomes $xf_{s+1}(x)$ and the right-hand side turns to be
\begin{align*}
&\sum\limits_{k=1}^n(-1)^{k-1}\binom{n}{k}\sum\limits_{l=0}^{\infty}\int_0^x\frac{\left(\frac{1}{2}\right)_l\cdots\left(\frac{1}{2}\right)_l\cdot(1-k)_l\cdot t^{2l}}{\left(\frac{3}{2}\right)_l\cdots\left(\frac{3}{2}\right)_l\cdot l!}dt\\
=&\sum\limits_{k=1}^n(-1)^{k-1}\binom{n}{k}\sum\limits_{l=0}^{\infty}\frac{\left(\frac{1}{2}\right)_l\cdots\left(\frac{1}{2}\right)_l\cdot(1-k)_l\cdot x^{2l+1}}{\left(\frac{3}{2}\right)_l\cdots\left(\frac{3}{2}\right)_l\cdot (2l+1)\cdot l!},
\end{align*}
which is
\begin{align*}
\sum\limits_{k=1}^n(-1)^{k-1}\binom{n}{k}\cdot x\cdot{_{s+2}F_{s+1}}\left(\left\{\frac{1}{2}\right\}^{s+1},1-k;\left\{\frac{3}{2}\right\}^{s+1};x^2\right).
\end{align*}
Hence, we prove that \eqref{3-1} holds for $s+1$.
\qed

Setting $x=1$ in Lemma \ref{lemma 3-1}, we get the evaluations of $\overline{H}_n(s)$. In the evaluation formulas, the values of the generalized hypergeometric functions at $1$ appear.

\begin{cor}
For $n,s\in \mathbb{N}$, we have
\begin{align*}
\overline{H}_n(s)=\sum\limits_{k=1}^n(-1)^{k-1}\binom{n}{k}{_{s+1}F_s}\left(\left\{\frac{1}{2}\right\}^{s},1-k;\left\{\frac{3}{2}\right\}^{s};1\right).
\end{align*}
\end{cor}

Note that if $s=1$, using the Chu-Vandermonde summation formula
$$_2F_1(-n,b;c;1)=\frac{(c-b)_n}{(c)_n},$$
we find that
$$\overline{H}_n(1)=\sum\limits_{k=1}^n(-2)^{k-1}\binom{n}{k}\frac{(k-1)!}{(2k-1)!!}.$$

Similarly as Lemma \ref{lemma 3-1}, we can prove the following result. And we omit the proof.

\begin{lem}\label{lemma 3-2}
For $n,s\in \mathbb{N}$, we have
\begin{align*}
\sum\limits_{k=0}^{n-1}\frac{(-1)^kx^{2k}}{(2k+1)^s}=\sum\limits_{k=1}^n(-1)^{k-1}\binom{n}{k}{_{s+1}F_s}\left(\left\{\frac{1}{2}\right\}^{s},1-k;\left\{\frac{3}{2}\right\}^{s};-x^2\right).
\end{align*}
\end{lem}

Let $x=1$ in Lemma \ref{lemma 3-2}, we get the evaluations of $\overline{H}_n(\overline{s})$. In the evaluation formulas, the values of the generalized hypergeometric functions at $-1$ appear.

\begin{cor}
For $n,s\in \mathbb{N}$, we have
\begin{align*}
\overline{H}_n(\overline{s})=\sum\limits_{k=1}^n(-1)^{k-1}\binom{n}{k}{_{s+1}F_s}\left(\left\{\frac{1}{2}\right\}^{s},1-k;\left\{\frac{3}{2}\right\}^{s};-1\right).
\end{align*}
\end{cor}

As a variant of the odd multiple harmonic sums, for $m,n\in\mathbb{N}$, we define
$$G_{m,n}:=\sum\limits_{k=0}^{n-1}\frac{1}{(2k+1)(2k+2)\cdots(2k+m)}.$$
Then we can represent $G_{m,n}$ by the special values of the hypergeometric function at $-1$.

\begin{prop}
For $m,n\in \mathbb{N}$, we have
\begin{align*}
G_{m,n}&=\sum\limits_{k=1}^n\frac{(-1)^{k-1}}{(m-1)!(m+k-1)}\binom{n}{k}{_2}F_1(1,1-k;m+k;-1).
\end{align*}
\end{prop}

\proof
Let $$h_m(x)=\sum\limits_{k=0}^{n-1}\frac{x^{2k+m}}{(2k+1)(2k+2)\cdots(2k+m)}.$$
Then we get
$$h^{(m)}(x)=\sum\limits_{k=0}^{n-1}x^{2k}=\sum\limits_{k=1}^n(-1)^{k-1}\binom{n}{k}(1-x^2)^{k-1}.$$
Integrating of both sides, we get
\begin{align*}
G_{m,n}&=\sum\limits_{k=1}^n(-1)^{k-1}\binom{n}{k}\int_0^1\int_0^{t_{m-1}}\cdots\int_0^{t_1}(1-x^2)^{k-1}dxdt_1\cdots dt_{m-1}\\
&=\sum\limits_{k=1}^n(-1)^{k-1}\binom{n}{k}\int_0^1\left(\,\int\limits_{x<t_1<\cdots<t_{m-1}<1}dt_1\cdots dt_{m-1}\right)(1-x^2)^{k-1}dx\\
&=\sum\limits_{k=1}^n\frac{(-1)^{k-1}}{(m-1)!}\binom{n}{k}\int_0^1(1-x)^{m+k-2}(1+x)^{k-1}dx\\
&=\sum\limits_{k=1}^n\frac{(-1)^{k-1}}{(m-1)!}\binom{n}{k}\sum\limits_{l=0}^{k-1}\binom{k-1}{l}\int_0^1x^l(1-x)^{m+k-2}dx\\
&=\sum\limits_{k=1}^n\frac{(-1)^{k-1}}{(m-1)!}\binom{n}{k}\sum\limits_{l=0}^{k-1}\binom{k-1}{l}B(l+1,m+k-1),
\end{align*}
where $B(a,b)$ denotes the beta function. Since $B(a,b)=\frac{(a-1)!(b-1)!}{(a+b-1)!}$ for positive integers $a,b$, we have
\begin{align*}
G_{m,n}&=\sum\limits_{k=1}^n\frac{(-1)^{k-1}}{(m-1)!}\binom{n}{k}\sum\limits_{l=0}^{k-1}\binom{k-1}{l}\frac{l!(m+k-2)!}{(m+k+l-1)!}\\
&=\sum\limits_{k=1}^n\frac{(-1)^{k-1}}{(m-1)!(m+k-1)}\binom{n}{k}\sum\limits_{l=0}^{k-1}\frac{(1)_l(1-k)_l(-1)^l}{l!(m+k)_l}\\
&=\sum\limits_{k=1}^n\frac{(-1)^{k-1}}{(m-1)!(m+k-1)}\binom{n}{k}{_2}F_1(1,1-k;m+k;-1).
\end{align*}
Hence, we complete the proof.
\qed

As applications, using the binomial inversion
$$f(n)=\sum\limits_{k=1}^n\binom{n}{k}g(k)\Longleftrightarrow g(n)=\sum\limits_{k=1}^n
(-1)^{n-k}\binom{n}{k}f(k),$$
we obtain the following summation formulas of the (generalized) hypergeometric functions.

\begin{cor}
For $m,n,s\in \mathbb{N}$, we have
\begin{align*}
&{_{s+1}}F_s\left(\left\{\frac{1}{2}\right\}^{s},1-n;\left\{\frac{3}{2}\right\}^{s};1\right)=\sum\limits_{k=1}^n(-1)^{k-1}\binom{n}{k}\overline{H}_k(s),\\
&{_{s+1}}F_s\left(\left\{\frac{1}{2}\right\}^{s},1-n;\left\{\frac{3}{2}\right\}^{s};-1\right)=\sum\limits_{k=1}^n(-1)^{k-1}\binom{n}{k}\overline{H}_k(\overline{s}),
\end{align*}
and
$${_2}F_1(1,1-n;m+n;-1)=(m-1)!(m+n-1)\sum\limits_{k=1}^n(-1)^{k-1}\binom{n}{k}G_{m,k}.$$
\end{cor}

\begin{rem}
Similarly, we can obtain evaluations of the (alternating) multiple harmonic sums of depth one.  For any $n,s\in\mathbb{N}$, we have
\begin{align*}
H_n(s)&=\sum\limits_{k=1}^n(-1)^{k-1}\binom{n}{k}{_{s+1}}F_s\left(\left\{1\right\}^{s},1-k;\left\{2\right\}^{s};1\right),\\
H_n(\overline{s})=&\sum\limits_{k=1}^n(-1)^{k-1}\binom{n}{k}{_{s+1}}F_s\left(\left\{1\right\}^{s},1-k;\left\{2\right\}^{s};-1\right),
\end{align*}
where $H_n(\overline{s})$ is the alternating harmonic sum defined by
$$H_n(\overline{s})=\sum\limits_{k=1}^n\frac{(-1)^{k-1}}{k^s}.$$
Using the binomial inversion formula, we may obtain the summation formulas
\begin{align*}
&{_{s+1}}F_s\left(\left\{1\right\}^{s},1-n;\left\{2\right\}^{s};1\right)=\sum\limits_{k=1}^n(-1)^{k-1}\binom{n}{k}H_k(s),\\
&{_{s+1}}F_s\left(\left\{1\right\}^{s},1-n;\left\{2\right\}^{s};-1\right)=\sum\limits_{k=1}^n(-1)^{k-1}\binom{n}{k}H_k(\overline{s}).
\end{align*}
\end{rem}

\end{document}